\documentclass{article}

\usepackage{graphicx,etoolbox, xcolor
}

\usepackage{hyperref}

\usepackage{hyperref}

\usepackage{amssymb}

\newcommand{\newc}{\newcommand}

\newc{\eqnoset}{\setcounter{equation}{0}}

\newcommand{\mref}[1]{(\ref{#1})}
\newcommand{\reflemm}[1]{Lemma~\ref{#1}}
\newcommand{\refrem}[1]{Remark~\ref{#1}}
\newcommand{\reftheo}[1]{Theorem~\ref{#1}}

\newcommand{\refcoro}[1]{Corollary~\ref{#1}}

\newcommand{\refsec}[1]{Section~\ref{#1}}

\newcommand{\beq}{\begin{equation}}
	\newcommand{\eeq}{\end{equation}}
\newcommand{\beqno}[1]{\begin{equation}\label{#1}}
	
	\newcommand{\barr}{\begin{array}}
		\newcommand{\earr}{\end{array}}
	
	\newc{\bearr}{\begin{eqnarray*}}
		\newc{\eearr}{\end{eqnarray*}}
	
	\newc{\bearrno}[1]{\begin{eqnarray}\label{#1}}
		\newc{\eearrno}{\end{eqnarray}}
	
	\newc{\non}{\nonumber}
	\newc{\nol}{\nonumber\nl}
	
	\newcommand{\bdes}{\begin{description}}
		\newcommand{\edes}{\end{description}}
	\newc{\benu}{\begin{enumerate}}
		\newc{\eenu}{\end{enumerate}}
	\newc{\btab}{\begin{tabular}}
		\newc{\etab}{\end{tabular}}

	\newtheorem{theorem}{Theorem}[section]
	\newtheorem{defi}[theorem]{Definition}
	\newtheorem{lemma}[theorem]{Lemma}
	\newtheorem{rem}[theorem]{Remark}
	\newtheorem{exam}[theorem]{Example}
	\newtheorem{propo}[theorem]{Proposition}
	\newtheorem{corol}[theorem]{Corollary}
	\newtheorem{conj}[theorem]{Conjecture}

	\newcommand{\btheo}[1]{\begin{theorem}\label{#1}}
		\newc{\brem}[1]{\begin{rem}\label{#1}\em}
			\newc{\bexam}[1]{\begin{exam}\label{#1}\em}
				\newc{\bdefi}[1]{\begin{defi}\label{#1}}
					\newcommand{\blemm}[1]{\begin{lemma}\label{#1}}
						\newcommand{\bprop}[1]{\begin{propo}\label{#1}}
							\newcommand{\bcoro}[1]{\begin{corol}\label{#1}}
								\newcommand{\btheoc}[1]{\begin{conj}\label{#1}}
									\newcommand{\etheo}{\end{theorem}}
								\newc{\etheoc}{\end{conj}}
							\newcommand{\elemm}{\end{lemma}}
						\newcommand{\eprop}{\end{propo}}
					\newcommand{\ecoro}{\end{corol}}
				\newc{\erem}{\end{rem}}
			\newc{\eexam}{\end{exam}}
		\newc{\edefi}{\end{defi}}
	
	\newc{\rmk}[1]{{\bf REMARK #1: }}
	\newc{\DN}[1]{{\bf DEFINITION #1: }}
	
	\newcommand{\bproof}{{\bf Proof:~~}}
	\newc{\eproof}{{\vrule height8pt width5pt depth0pt}\vspace{3mm}}
	
	\newc{\bfrac}[2]{\dspl{\frac{#1}{#2}}}

	\newc{\nid}{\noindent}
	
	\newcommand{\dspl}{\displaystyle}
	\newc{\grad}{\nabla}
	\newc{\Div}{\mbox{div}}
	\newc{\pdt}[1]{\dspl{\frac{\partial{#1}}{\partial t}}}
	\newc{\pdn}[1]{\dspl{\frac{\partial{#1}}{\partial \nu}}}
	\newc{\pdNi}[1]{\dspl{\frac{\partial{#1}}{\partial \mathcal{N}_i}}}
	\newc{\pD}[2]{\dspl{\frac{\partial{#1}}{\partial #2}}}
	\newc{\dt}{\dspl{\frac{d}{dt}}}
	\newc{\bdry}[1]{\mbox{$\partial #1$}}
	\newc{\sgn}{\mbox{sign}}
	
	\newc{\Hess}[1]{\frac{\partial^2 #1}{\pdh z_i \pdh z_j}}
	\newc{\hess}[1]{\partial^2 #1/\pdh z_i \pdh z_j}

	\newc{\ag}{\alpha}
	\newc{\bg}{\beta}
	\newc{\cg}{\gamma}\newc{\Cg}{\Gamma}
	\newc{\dg}{\delta}\newc{\Dg}{\Delta}
	\newc{\eg}{\varepsilon}
	\newc{\zg}{\zeta}
	\newc{\thg}{\theta}
	\newc{\llg}{\lambda}\newc{\LLg}{\Lambda}
	\newc{\kg}{\kappa}
	\newc{\rg}{\rho}
	\newc{\sg}{\sigma}\newc{\Sg}{\Sigma}
	\newc{\tg}{\tau}
	\newc{\fg}{\phi}\newc{\Fg}{\Phi}
	\newc{\vfg}{\varphi}
	\newc{\og}{\omega}\newc{\Og}{\Omega}
	\newc{\pdh}{\partial}
	
	\newc{\ccG}{{\cal G}}

	\newc{\ii}[1]{\int_{#1}}
	\newc{\iidx}[2]{{\dspl\int_{#1}~#2~dx}}
	\newc{\bii}[1]{{\dspl \ii{#1} }}
	\newc{\biii}[2]{{\dspl \iii{#1}{#2} }}
	\newc{\su}[2]{\sum_{#1}^{#2}}
	\newc{\bsu}[2]{{\dspl \su{#1}{#2} }}

	\newc{\biiom}[1]{{\dspl\int_{\bdrom}~ #1 ~d\sg}}
	\newc{\io}[1]{{\dspl\int_{\Og}~ #1 ~dx}}
	\newc{\bio}[1]{{\dspl\int_{\bdrom}~ #1 ~d\sg}}
	\newc{\bsir}{\bsu{i=1}{r}}
	\newc{\bsim}{\bsu{i=1}{m}}
	
	\newc{\iibr}[2]{\iidx{\bprw{#1}}{#2}}
	\newc{\Intbr}[1]{\iibr{R}{#1}}
	\newc{\intbr}[1]{\iibr{\rg}{#1}}
	\newc{\intt}[3]{\int_{#1}^{#2}\int_\Og~#3~dxdt}
	
	\newc{\itQ}[2]{\dspl{\int\hspace{-2.5mm}\int_{#1}~#2~dz}}
	\newc{\mitQ}[2]{\dspl{\rule[1mm]{4mm}{.3mm}\hspace{-5.3mm}\int\hspace{-2.5mm}\int_{#1}~#2~dz}}
	\newc{\mitQQ}[3]{\dspl{\rule[1mm]{4mm}{.3mm}\hspace{-5.3mm}\int\hspace{-2.5mm}\int_{#1}~#2~#3}}
	
	\newc{\mitx}[2]{\dspl{\rule[1mm]{3mm}{.3mm}\hspace{-4mm}\int_{#1}~#2~dx}}
	\newc{\mitmu}[2]{\dspl{\rule[1mm]{3mm}{.3mm}\hspace{-4mm}\int_{#1}~#2~d\mu}}
	\newc{\iidmu}[2]{\iidx{#1}{#2}}
	
	\newc{\iidm}[3]{{\dspl\int_{#1}~#2~d #3}}
	
	\newc{\itQmu}[2]{\dspl{\int\hspace{-2.5mm}\int_{#1}~#2~d\mu}}
	\newc{\mitQmu}[2]{\dspl{\rule[1mm]{4mm}{.3mm}\hspace{-5.3mm}\int\hspace{-2.5mm}\int_{#1}~#2~d\mu}}

	\newc{\mitQq}[2]{\dspl{\rule[1mm]{4mm}{.3mm}\hspace{-5.3mm}\int\hspace{-2.5mm}\int_{#1}~#2~d\bar{z}}}
	\newc{\itQq}[2]{\dspl{\int\hspace{-2.5mm}\int_{#1}~#2~d\bar{z}}}

	\newc{\pder}[2]{\dspl{\frac{\partial #1}{\partial #2}}}

	\newc{\bdrom}{\bdry{\Og}}

	\newc{\bilhom}{\mbox{Bil}(\mbox{Hom}(\RR^{nm},\RR^{nm}))}
	\newc{\VV}[1]{{V(Q_{#1})}}
	
	\newc{\ccA}{{\mathcal A}}
	\newc{\ccB}{{\mathcal B}}
	\newc{\ccC}{{\mathcal C}}
	\newc{\ccD}{{\mathcal D}}
	\newc{\ccE}{{\mathcal E}}
	\newc{\ccH}{\mathcal{H}}
	\newc{\ccF}{\mathcal{F}}
	\newc{\ccI}{{\mathcal I}}
	\newc{\ccJ}{{\mathcal J}}
	\newc{\ccK}{{\mathcal K}}
	\newc{\ccP}{{\mathcal P}}
	\newc{\ccQ}{{\mathcal Q}}
	\newc{\ccR}{{\mathcal R}}
	\newc{\ccS}{{\mathcal S}}
	\newc{\ccT}{{\mathcal T}}
	\newc{\ccX}{{\mathcal X}}
	\newc{\ccY}{{\mathcal Y}}
	\newc{\ccZ}{{\mathcal Z}}
	
	\newc{\bb}[1]{{\mathbf #1}}
	
	\newc{\myprod}[1]{\langle #1 \rangle}
	\newc{\mypar}[1]{\left( #1 \right)}
	
	\newc{\BLLg}{\mathbf{\LLg}}
	
	\newc{\mA}{\mathbf{A}}
	\newc{\mB}{\mathbf{B}}
	\newc{\mC}{\mathbf{C}}
	\newc{\mD}{\mathbf{D}}
	\newc{\mE}{\mathbf{E}}
	\newc{\mF}{\mathbf{F}}
	\newc{\mJ}{\mathbf{J}}
	\newc{\mG}{\mathbf{G}}
	\newc{\mP}{\mathbf{P}}
	\newc{\mR}{\mathbf{R}}
	\newc{\mQ}{\mathbf{Q}}
	\newc{\mX}{\mathbf{X}}
	\newc{\muu}{\mathbf{u}}
	\newc{\mvv}{\mathbf{v}}
	
	\newc{\mllg}{\mathbb{\lambda}}
	\newc{\mLLg}{\mathbf{\LLg}}

	\newc{\lspn}[2]{\mbox{$\| #1\|_{\Lsp{#2}}$}}
	\newc{\Lpn}[2]{\mbox{$\| #1\|_{#2}$}}
	\newc{\Hn}[1]{\mbox{$\| #1\|_{H^1(\Og)}$}}

	\newc{\mynorm}[2]{\| #1\|_{#2}}

	\newcommand{\RR}{{\rm I\kern -1.6pt{\rm R}}}

	\newc{\itQQ}[2]{\dspl{\int_{#1}#2\,dz}}
	\newc{\mmitQQ}[2]{\dspl{\rule[1mm]{4mm}{.3mm}\hspace{-4.3mm}\int_{#1}~#2~dz}}
	\newc{\MmitQQ}[2]{\dspl{\rule[1mm]{4mm}{.3mm}\hspace{-4.3mm}\int_{#1}~#2~d\mu}}

	\newc{\MUmitQQ}[3]{\dspl{\rule[1mm]{4mm}{.3mm}\hspace{-4.3mm}\int_{#1}~#2~d#3}}
	\newc{\MUitQQ}[3]{\dspl{\int_{#1}~#2~d#3}}

	\newc{\mccP}{\mathbb{P}}
	\newc{\mccK}{\mathbb{K}}
	
	\newc{\DKTmU}{\mccK(U)}
	\newc{\DKTmUold}{(K_U(U)^{-1})^T}
	
	\newc{\myPi}{\mathbf{W}}
	\newc{\myIbar}{\bar{\ccI}_1}
	\newc{\myIhat}{\hat{\ccI}_1}
	\newc{\myIbreve}{\breve{\ccI}_0}
	
	\newc{\mmk}{\mathbf{k}}

	\newcommand{\mg}{\mathbf{g}}

	\newc{\mfu}{\mathbf{f_u}}
	\newc{\mh}{\mathbf{h}}
	\newc{\mb}{\mathbf{b}}
	\newc{\mf}{\mathbf{f}}

	\newcommand{\barrl}[2]{\barr{ll}\lefteqn{#1}\hspace{#2}&\\}

	\newc{\twomatrix}[1]{\left[\barr{cc}#1\earr\right]}
	\newc{\threematrix}[1]{\left[\barr{ccc}#1\earr\right]}

	\newc{\mN}{\mathbf{N}}
	\newc{\mI}{\mathbf{I}}
	\newc{\mH}{\mathbf{H}}

	\newc{\mk}{\mathbf{k}}
	\newc{\mr}{\mathbf{r}}

	\newc{\DIAGM}[2]{\left[\barr{ccc}#1&0\ldots&0\\
		\vdots&\ddots&\vdots\\0&\ldots0&#2\earr \right]}
	\newc{\DiagM}[2]{\mbox{diag}\left[#1
		\cdots #2 \right]}
	\newc{\vVEC}[2]{\left[\barr{c}#1\\
		\vdots\\#2\earr \right]}
	\newc{\hVEC}[2]{\left[#1
		\cdots #2 \right]}
	
	\newc{\mq}{\mathbf{q}}

	\newc{\msys}[1]{\left\{\barr{l}#1\earr
		\right.}
	\newc{\msysa}[1]{\left\{\barr{ll}#1\earr
		\right.}
	
	\newc{\bbM}{\mathbb{M}}
	\newc{\mat}[1]{\left[\barr{cc}#1\earr\right]}
	
	\newc{\me}{\mathbf{e}}

	\newc{\vecc}[2]{\left[\barr{cc}#1\\#2\earr\right]}
	\newc{\mL}{\mathbb{L}}
	
	\newc{\cO}{{\cal O}}
	
	\newc{\cM}{{\cal M}}
	
	\newc{\myega }{\eg_0(R)}
	\newc{\myeg}{\eg_1(\eg_*)}
	\newc{\myegp}{\hat{\eg}_1(\eg_*)}

	\newc{\diagA}{\mathbb{A}_d}
	\newc{\mBB}{\mathbb{B}}
	\newc{\MLT}[1]{{\cal M}_{lt}(\Og,#1)}
	\newc{\ALT}[1]{{\cal A}_{l}(\Og,#1)}
	\newc{\mM}{\mathbb{M}}
	
	\newc{\diag}[1]{\mbox{diag}(#1)}
	\newc{\off}[1]{\mbox{offdiag}(#1)}
	\newc{\mT}{\mathbb{T}}

	\usepackage{amssymb}

	\newc{\idmu}[2]{{\dspl\int_{#1}~#2~d\mu}}
	\newc{\idllg}[2]{{\dspl\int_{#1}~#2~d\llg}}
	
	\newc{\wsp}[2]{\dspl{\int}_{#1}\dspl{\int}_{#1}\frac{|#2(x)-#2(y)|^p}{|x-y|^{N+sp}}dxdy}
	\newc{\Wsp}[2]{\|#2\|_{W^{s,p}(#1)}}
	
	\newc{\dsg}{(-\Delta)^{\frac{\sg}{2}}}
	\newc{\dgeng}[1]{(-\Delta)^{\frac{#1}{2}}}
	
	\newc{\dcg}{(-\Delta)^{\frac{\cg}{2}}}
	
	\newc{\ds}{(-\Delta)^{\frac{s}{2}}}
	
	\newc{\dsab}{(-\Delta_{||})^{\frac{s}{2}}}
	\newc{\dsgab}{(-\Delta_{||})^{\frac{\sg}{2}}}

	\newc{\midx}[2]{\dspl{\rule[1mm]{3mm}{.3mm}\hspace{-4mm}\int_{#1}~#2~dx}}
	
	\newc{\midz}[3]{\dspl{\rule[1mm]{3mm}{.3mm}\hspace{-4mm}\int_{#1}~#2~ #3}}
	
	\newc{\idz}[3]{\dspl{\int_{#1}~#2~ #3}}
	
	\newc{\cfs}{C_{FS}}
	\newc{\cps}{C_{PS}}
	\newc{\nc}{\mathbf{c}}

	\newc{\xoa}{{\bf Delete?}}
	\newc{\coilai}{{\bf ???}}

\newc{\myrefc}[1]{{\bf refer #1}}
\newc{\mybook}[1]{#1}

\begin{document}

	\vspace*{-.8in}
	\begin{center} {\LARGE\em Gagliardo-Nirenberg type inequalities with a BMO term and fractional Laplacians}
		
	\end{center}

	\vspace{.1in}
	
	\begin{center}

		{\sc Dung Le}{\footnote {Department of Mathematics, University of
				Texas at San
				Antonio, One UTSA Circle, San Antonio, TX 78249. {\tt Email: Dung.Le@utsa.edu}
				}}

	\end{center}

	\begin{abstract} An improvement of  a {\em Global (strong) Gagliardo-Nienberg inequality with a BMO term} is established by replacing local derivatives by {\em  and fractional Laplacians.} Local versions are also given.
		
	 \end{abstract}

\section{Introduction} \label{fracintro} In this paper we provide some improvements of  the global (strong) Gagliardo-Nienberg inequality with a BMO term in \cite{SR} and
its local (weak) Gagliardo-Nienberg inequality with a BMO term in \cite{SR,dlebook1} by extending these results in using the non-local operators fractional Laplacians.

We will show that: Let $H,u$ be  such that $u\in BMO(B_R)$, $B_R$ denotes a ball in $\RR^N$ with radius $R>0$, and  $u, \dsg u\in L^2(\RR^N),\; |H|^{2p} \in L^1(\RR^N)$, $0<\sg\le 1, p\ge1$. Then,
$$\iidx{\RR^N}{|H|^{2p}|\dsg u|^2}\le C(N,p,\sg)\|u\|_{BMO(B_R)}^2\iidx{\RR^N}{|H|^{2p-2}|\dsg H|^2}.$$

The above inequality generalizes the global (strong) Gagliardo-Nienberg inequality with a BMO term in \cite{SR}, which  proves that $\|Df\|_{L^{2p}(\RR^N)}^2\le C\|f\|_{BMO}\|f\|_{W^{2,p}(\RR^N)}$. It is also more versatile than the strong version as it allows two functions $H,u$ and uses the non-local operators fractional Laplacians $\dsg$ instead of the local operator $D$. In fact, by choosing $H,u$ appropriately, we can derive the strong version from this inequality. 

It is possible to incorporate fractional diffusion into nonlinear diffusion models (see \cite{vasquez}). The regularity of weak solutions becomes even harder (\cite{zua}). We are able to extend the local (weak) Gagliardo-Nienberg inequality with a BMO term to deal with fractional Laplacians $\dsg$ and this can combine with the above techniques to show similar H\"older regularity results. Indeed, because the difference operator $\dg$ commutes with $\dsg$ so that the previous arguments (e.g. \cite{dlebook1}) can be applied. However, the power rules (for the local derivative $D$ in \cite{dlebook1})  are problematic for fractional Laplacians. Nevertheless, the connection between a local property (BMO) and a non local operators (fractional Laplacians) is of interest by itself and can be useful for PDE's analysis as in \cite{vasquez}.

In \refsec{Wspsapce}, we gather some basic information on the fractional Sobolev spaces and fractional Laplacians $\dsg$. A good reference of these can be found in \cite{NPV}.  We also collect our main tools in Harmonic analysis used in our proofs. An excellent source for these is the book \cite{Grafa}.

\refsec{newWGNsec} presents our main result: the (local weak) Gagliardo-Nienberg inequality for fractional Laplacians $\dsg$ with a BMO term-and its proof. Here, we use the famous Fefferman-Stein inequality as in \cite{SR} and nontrivial modifications.

To the best of our knowledge, there is no theory available for strong solutions to diffusion equations using fractional Laplacians and thus a {\em local weak Gagliardo-Nienberg inequality with a BMO term and fractional Laplacians $\dsg$}  could be a useful tool to study the regularity of weak solutions to such euqations/systems.  For the readers, who are not familiar with the difficulties of (nonlocal) fractional Laplacians,  they can easily use the standard (local) derivative $D$ instead. The main ideas of the proofs are the same.

\section{Preliminaries}\label{Wspsapce}\eqnoset

\subsection{$W^{s,p}$ spaces and frational Laplacians}
As usual, we define the $W^{s,p}(\Og)$ spaces (when $s\in(0,1)$ and $\Og$ is an extension domain, e.g.  $\Og=B_R$) with finite Gagliardo seminorm $$\Wsp{\Og}{u}=\left(\wsp{\Og}{u}\right)^\frac{1}{p}.$$
We have the  well known inequality for $p>1$, $s\in(0,1)$ and $sp<N$ 
$$\left( \midz{\Og}{|u-u_{\Og}|^{\frac{Np}{N-sp}}}{dy} \right)^\frac{N-sp}{Np}\le C|\Og|^\frac{s}{N} \left(\frac{1}{|\Og|}\wsp{\Og}{u}\right)^\frac{1}{p}.$$

The fractional Laplacian is usually defined (pointwise) by \beqno{fraclappv}\ds u(x)= C_{N,s}\mbox{P.V.}\idz{\RR^N}{\frac{u(x)-u(y)}{|x-y|^{N+s}}}{dy}= C_{N,s}\lim_{r\to0}\idz{\{|x-y|>r\}}{\frac{u(x)-u(y)}{|x-y|^{N+s}}}{dy}\eeq
where $C_{N,s}=2^{s-1}s\Gamma((N+s)/1)/(\pi^{N/2}\Gamma(1-s/2))$. 
\brem{prodineq} The pointwise definition \mref{fraclappv} explains why $|u(y)| \le C(\Og)|u(x)|$ for compact supp(u). Thus, we have power rule {\em inequality} because $|u(y)|\le C(B_x)|u(x)|$ locally. $C(B_x)$ may depend on $K= \|u\|_W^{s,p}$ (or some other norms, we can take $K=1$ by scalings).
Similarly, we have product rule {\em inequality} for $\dsg(uv)$ if $supp(u)$ is compact and $v \in  C^1(R^N)$.
\erem

We also recall for $u\in \mathfrak{S}$ (the set of tempered distributions on $\RR^N$ and using Fourier's transform $\mathfrak{F}$, e.g.\cite{NPV}) the fractional Laplacian is also defined by 
\beqno{fraclapF}(-\Delta)^s u=\mathfrak{F}^{-1}(|\zeta|^{2s}\mathfrak{F}u) \mbox{ so } (-\Delta)^\frac{s}{2} u=\mathfrak{F}^{-1}(|\zeta|^{s}\mathfrak{F}u).\eeq

We also recall ($\dsg$ is defined by the two equivalent (in $L^2$ norm) \mref{fraclappv} and \mref{fraclapF})
$$H^\sg(\RR^N)=\{u:\, \|u\|_{L^2(\RR^N)}+\|\dsg u\|_{L^2(\RR^N)}<\infty\},$$
and $H^\sg(\RR^N)=W^{\sg,2}(\RR^N)$.

\brem{normequiv} By \cite{MK} the norms of $H^\sg(\RR^N)=W^{\sg,2}(\RR^N)$ are equivalent with $\dsg$ given by \mref{fraclappv} and \mref{fraclapF}. This and the parallelogram formula for product spaces  imply
that the {\em absolutes} of $\myprod{u,v}=\myprod{u,v}_{L^2}+\myprod{\dsg u, \dsg v}_{L^2}$ are also equivalent. \erem

When $p=2$, the  Poincar\'e-Sobolev's inequality is valid for the standard fractional derivatives using Fourier's transform and extension \cite[Sect. 3.1]{NPV}.

\blemm{HsSobo} If $\Og$ is an extension domain then for $s\in(0,1]$
$$\left( \midz{\Og}{|u-u_{\Og}|^{\frac{2N}{N-2s}}}{dy} \right)^\frac{N-2s}{2N}\le  \cps(N,s)|\Og|^\frac{s}{N}\left(\midz{\Og}{|\ds u|^2}{dy}\right)^\frac{1}{2}. $$
\elemm

Thanks to this and \refrem{normequiv}, we see that the absolute of $\myprod{\dsg u, \dsg v}_{L^2}$ (with $\dsg$ given by \mref{fraclappv} and \mref{fraclapF}) are also equivalent for compact support $u,v$.
	
We also use the fact that in $H^s(\RR^N)=W^{s,2}(\RR^N)$ the integration by parts is available for $\dsg$ (defined by \mref{fraclapF}). As we cannot find the reference for this following  elementary result, we present here the simple proof.

\blemm{lapintlem}
If $ \ag, \bg\in[0,\infty)$ and $f,u$ are tempered distributions, then 
$$\iidx{\RR^N}{\myprod{ f, \overline{ \dgeng{\ag+\bg} u} }}=\iidx{\RR^N}{\myprod{\dgeng{\ag} f, \overline{\dgeng{\bg} u} }}.$$

\elemm

This also shows that $\dsg = \dgeng{\sg-\bg}\dgeng{\bg}$ for $\sg>\bg>0$.

\bproof Indeed, if we use the Fourier transform to define $\dsg$ and Parseval's theorem twice then 

$$\barr{lll}\iidx{\RR^N}{\myprod{\dgeng{\ag}f ,\overline{\dgeng{\bg}u}} }&=& \iidx{\RR^N}{\myprod{\mathfrak{F}^{-1}(|\zeta|^\ag \mathfrak{F}(f)),\overline{ \mathfrak{F}^{-1}(|\zeta|^\bg \mathfrak{F}(u)) }}}\\
\mbox{(Parseval)}&=& \idz{\RR^N}{\myprod{|\zeta|^\ag \mathfrak{F}(f), \overline{|\zeta|^\bg \mathfrak{F}(u)}}{d\zeta} } = \idz{\RR^N}{ \myprod{\mathfrak{F}(f), \overline{|\zeta|^{\ag+\bg} \mathfrak{F}(u)} }}{d\zeta}\\
\mbox{(Parseval)}&=& \iidx{\RR^N}{ \myprod{\mathfrak{F}^{-1}(\mathfrak{F}(f)), \overline{\mathfrak{F}^{-1}(|\zeta|^{\ag+\bg}\mathfrak{F}(u))} }}\\
&=&\iidx{\RR^N}{\myprod{f, \overline{\dgeng{\ag+\bg} u} }}.
\earr
$$
The proof is complete. \eproof

\subsection{Notations and tools from Harmonic analysis}\label{mainres}

For any measurable subset $A$ of $\Og$  and any  locally integrable function $U:\Og\to\RR^m$ we denote by  $|A|$ the measure of $A$ and $U_A$ the average of $U$ over $A$. That is, $$U_A=\mitx{A}{U(x)} =\frac{1}{|A|}\iidx{A}{U(x)}.$$

We now recall some well known notions from Harmonic Analysis.

A function $f\in L^1(\Og)$ is said to be in $BMO(\Og)$ if \beqno{bmodef} [f]_{*}:=\sup_Q\mitx{Q}{|f-f_Q|}<\infty.\eeq We then define $$\|f\|_{BMO(\Og)}:=[f]_{*}+\|f\|_{L^1(\Og)}.$$

We also recall the Hardy-Littlewood maximal function
$$M(f)(x)=Mf(x):=\sup_{r>0}\midz{B_r(x)}{|f(y)|}{dy}$$
and the Hardy-Littlewood inequality: if $q>1$ and $w$ belongs to the class $A_\cg$
\beqno{HL}\iidx{\Og}{|M(f)(x)|^qw(x)}\le C(q)\iidx{\Og}{|f(x)|^qw(x)}.\eeq

For $\cg\in(1,\infty)$ we say that a nonnegative locally integrable function $w$ belongs to the class $A_\cg$ or $w$ is an $A_\cg$ weight if the quantity
$$ [w]_{\cg} := \sup_{B\subset\Og} \left(\mitx{B}{w}\right) \left(\mitx{B}{w^{1-\cg'}}\right)^{\cg-1} \quad\mbox{is finite}.$$
Here, $\cg'=\cg/(\cg-1)$ and the supremum is taken over all cubes $B$ in $\Og$.

Here and throughout this chapter, we write $B_R(x)$ for a cube centered at $x$ with side length $R$ and sides parallel to the standard axes of $\RR^N$. We will omit $x$ in the notation $B_R(x)$ if no ambiguity can arise. We denote by $l(B)$ the side length (radius) of a cube (ball) $B$ and by $\tau B$ the cube (ball) which is concentric with $B$ and has side length (radius) $\tau l(B)$.

One of the key ingredients of our proof is the duality between the Hardy space $\ccH^1(\RR^N)$ and $BMO(\RR^N)$ space. This is the famous Fefferman-Stein theorem (see \cite{fst}). It is well known that the norm of the Hardy space can be defined by
$$\|g\|_{\ccH^1(\RR^N)}=\|g\|_{L^1(\RR^N)}+\|M_*g\|_{L^1(\RR^N)}.$$ Here, $$M_{*}g =\sup_{\eg>0}|\fg_\eg*g|$$ where  $\fg_\eg=\eg^{-N}\fg(\frac{x}{\eg})$ with $\fg\in C^\infty(\RR^N)$ which has support in the ball $B_1(y)$ centered at $y$ and $\int_{\RR^N}\fg(y)dy=1$. The definition does not depend on the choice of $\fg$ (\cite{fst}) so that throughout this paper, by the properties of $\fg$ and scaling, for a given $\sg\in(0,1]$ we will always assume that \beqno{fgprop}0\le \fg_\eg(x)\le \eg^{-N},\; |\dsg \fg_\eg(x)|\le \eg^{-N-\sg} \mbox{ for all $x\in \RR^N$}.\eeq
The last property easily comes from the diameter of $\mbox{supp}\fg_\eg$ and the definition \mref{fraclappv}  and a change of variables.

The famous Fefferman-Stein theorem states that there is a constant $\cfs(N)$ such that \beqno{FSineq}\left|\iidx{\RR^N}{f(x)g(x)}\right|\le \cfs(N)\|f\|_{BMO}\|g\|_{{\cal H}^1(\RR^N)}.\eeq

We will also use the definition of the  {\em centered}  Hardy-Littlewood maximal operator acting on functions $F\in L^1_{loc}(\Og)$
\beqno{maximal} M(F)(y) = \sup_\eg\{\mitx{B_\eg(y)}{F(x)}\,:\, \eg>0 \mbox{ and } B_\eg(y)\subset\Og\}.\eeq



\section{A new local weak Gagliardo-Nirenberg inequality} \label{newWGNsec} \eqnoset Let $\og$ be a function with compact support in $\RR^N$. We  assume that $\dsg H,\dsg u, \dsg\og$ are  sufficiently integrable.  The new inequality presented here is another version of  the weak GNBMO inequality in \cite{dlebook1} when $\sg=1$ and also allows us to consider weak solutions and study their local regularities. The proof is also quite similar but needs some major modifications.  We also discuss several  variations of it which may be of interest in other applications.
We denote 
\newc{\mmyIbreve}{\mathbf{\myIbreve}}

\beqno{Idefz} \mI_1:=\iidx{\RR^N}{|H|^{2p}|\dsg u|^2\og^2},\;
\mI_2:=\iidx{\RR^N}{|H|^{2p-2}|\dsg H|^2\og^2},\eeq
\beqno{Idef2z}\mmyIbreve:=\iidx{\Og}{|H|^{2p}}.\eeq 

\newc{\Ofam}{{\mathfrak{O}_\og}}

{\bf Definition of the family $\Ofam$:} we partition $\Og$ into finitely many {\em sets } $B^i$'s, which has a finite intersection property, and choose $\og_i$'s such that $\og^2\le\sum_i\og^2_i$,  $\mbox{supp}\og_i\subset B^i$, and  $|\dsg \og_i|\le \LLg_i$ for some finite numbers $\LLg_i$'s.

\btheo{GNlobalz}  
Suppose that $p\ge1$ and  $\myega, \mI_1,\mI_2,\mmyIbreve$ are finite. \xoa For any given   $\eg_*>0$ we have a contant $C(N,p,\sg)$ such that
\beqno{GNglobalestz}\mI_1\le C(N,p,\sg)(\|u\|_{BMO(\Og_R)}+\eg_*)^2 \mI_2+E_\Ofam(u)\mmyIbreve+\iidx{\RR^N}{|H|^{2p}|\dsg\og|^2},\eeq
where \beqno{E*def}E_\Ofam(u)=\sup_{B^i\in\Ofam}|B^i|^{-1}\iidx{B^i}{[|\dsg u|^2+|u|^2\LLg_i^2]}.\eeq
\etheo

The estimate \mref{GNglobalestz} obviously yields the weak inequality in \cite{dlebook1} when $\sg=1$ and $B^i$'s are simply balls. The quantity $E_\Ofam(u)$ depends on the geometry of $B^i$'s of $\Ofam$ which will be investigated further (to be more precise $E_\Ofam \sim |B^i|^{-1}\|\dsg (u\og_i)\|_{L^2(\RR^N)}$ ). In fact, $\og$  in our local weak inequality of \reftheo{GNlobalz} can be removed, keeping only $\mI_1, \mI_2$, if $\Og$ is an extended domain (as $B_R$ (see \refcoro{locwGNcoro})). However, due to the common nature of the local analysis of the regularity of weak solutions to PDE's, one has $\og$ in the consideration and \reftheo{GNlobalz} is sufficient.

The proof consists some non trivial uses of the Fefferman-Stein theorem and the Sobolev inequality in \reflemm{HsSobo}.

\subsection{The proof:}

We prepare the proof by introducing the following functions $\mg_1,\mg_2,\mg_3$ with $\mg_i=\dsg (U_i)$, where (setting  $h=|H|^{p-1}H$) and for any $y\in \RR^N$ and $\eg>0$
$$U_1= h\dsg u\left(h-\mitx{B_{\eg}(y)}{h}\right),\;U_2= (h-\mitx{B_{\eg}(y)}{h})\dsg u\mitx{B_{\eg}}{h},$$
$$U_3=\dsg u\left(\mitx{B_{\eg}(y)}{h}\right)^2.$$

We see that $|H|^{2p}\dsg u=h^2\dsg u=U_1+U_2+U_3$ so that by integration by parts,  writing $|\dsg u|^2=\myprod{\dsg u,\dsg u}$ and \reflemm{lapintlem}, and product rule {\em inequality} (\refrem{prodineq})  \beqno{keyGNwsp}\mI_1=\iidx{\RR^N}{\dsg(|H|^{2p}\dsg u \og^2)u}\le \mathfrak{J}_1+\mathfrak{J}_2,\eeq where
\beqno{JJdef}\mathfrak{J}_1:=\left|\iidx{\RR^N}{(\mg_1+\mg_2) u\og^2}\right|,\; \mathfrak{J}_2:=\left|\iidx{\RR^N}{\mg_3 u\og^2}\right|+
\iidx{\RR^N}{|H|^{2p}|\dsg u||\og||\dsg\og|}. \eeq

Thus, $\mI_1$ can be  estimated by $\mathfrak{J}_i$'s. To estimate the  terms $\mathfrak{J}_i$'s in \mref{keyGNwsp}, we will use Fefferman-Stein theorem and  have the following lemmas.

{\em We emphasize that in the proof below we will use the Fefferman-Stein theorem to estimate $\mathfrak{J}_1$ and we will take the supremum in $\eg>0$ to do so. Meanwhile, we will not do so  for $h_{B_\eg(y)}$ in $\mg_3$ (in $\mathfrak{J}_2$). Thus, we can take any $\eg>0$ in estimating $\mathfrak{J}_2$.
}

\brem{defs-interchange} By \refrem{normequiv}  ,  the absolute of inner products $\myprod{\dsg u, \dsg v}$ are equivalent when $\dsg$ are defined by \mref{fraclapF} or \mref{fraclappv}. Therefore 2 definitions can be used interchangeably later as we will need integration by parts formula (\reflemm{lapintlem}) and product rule {\em inequality} (\refrem{prodineq}).
\erem
First of all, by the Fefferman-Stein theorem again ($\mbox{supp}\og\subset\Og_R$), $$\left|\iidx{\RR^N}{(\mg_1+\mg_2) u\og^2}\right|\le \cfs(N)\|u\|_{BMO{(B_R)}}\|(\mg_1+\mg_2)|\|_{\ccH^1(\Og_R)}.$$ 
We estimate $\|(\mg_1+\mg_2)|\|_{\ccH^1(\Og_R)}$ in the following lemmas.

\blemm{g12lemM} We have \beqno{newg1estzz}\iidx{\Og_R}{\sup_{\eg>0}|(\mg_1+\mg_2)*\fg_\eg|\og^2} \le \cps(N,\sg) \mI_1^\frac12\left\|M\left(|H|^{2(p-1)}|\dsg H|^{2}\right)\right\|_{L^1(\Og_R)}^\frac12.\eeq
\elemm

\bproof Let us consider $\mg_1$ first and estimate the term $\iidx{\Og}{\sup_\eg|\mg_1*\fg_\eg|}$. From \mref{fgprop},  we have $|\dsg \fg_\eg|\le \eg^{-N-\sg}$.
For any $y\in\Og$, we use integration by parts (\mref{fraclapF} via \reflemm{lapintlem}), ignoring $\og^2$ for simplicity  because the integrals involving $\og, \dsg \og$ in \mref{keyGNwsp} can be treated the same way, see  \reflemm{g3lem} below. Partition again $\Og$ into finitely many {\em balls } $B^\eg$'s, which has a finite intersection property and $\mbox{supp}\fg_\eg\subset B_\eg$. Using the definition \mref{fraclappv} and property of $\fg_\eg$, \mref{fgprop}, and then H\"older's inequality for any $s>1$, we have the following
$$\barr{lll}|\mg_1*\fg_\eg(y)|&=&\left|\iidx{B_\eg(y)}{\dsg \fg_\eg(x-y)(h-h_{B_\eg(y)})h\dsg u  }\right|\\&\le&\frac{C_1}{\eg^\sg}
\left|\mitx{B_\eg(y)}{|h-h_{B_\eg(y)}||h\dsg u|}\right|
\\&\le&\frac{C_1}{\eg^\sg}
\left(\mitx{B_\eg(y)}{|h-h_{B_\eg(y)}|^s}\right)^\frac1s 
\left(\mitx{B_\eg(y)}{|h\dsg u|^{s'}}\right)^\frac1{s'}.
\earr$$

Since $N\ge2$ and $\sg\in(0,1)$, we can take  $s=\frac{2N}{N-2\sg}$
 then $s_*=\frac{Ns}{N+\sg s}=2$. As $s_*=2$ and $\Og$ is an extension domain we can use Sobolev-Poincar\'e's inequality in \reflemm{HsSobo} here. 
We have
\beqno{hestz}\frac{C_1}{\eg^\sg}
\left(\mitx{B_\eg}{|h-h_{B_\eg}|^s}\right)^\frac1s \le C\cps(N,\sg)\Psi_2^\frac{1}{2},\eeq where $\Psi_2(y)=M(|H|^{2(p-1)}|\dsg H|^{2})(y)$. Here,  as $s_*=2$, we use the fact that 
$$|\dgeng{\sg} h| \le c(p)|H|^{p-1}|\dgeng{\sg}H|, \quad h=|H|^{p-1}H.$$ This can be seen easily by using the formula (see also \cite{Garo}) $$f(u)=\int_0^1f'(su)uds$$ and the Fourier transform on $W^{\sg,2}(\RR^N)$ (see \cite{NPV}), Fubini's theorem and an approximation of $f(u)=|u|^{p-1}u$ by piece-wise linear  functions. We omit elementary and technical details.

Setting $\Psi_3(y)=M(|h\dsg u|^{s'})(y)$ and putting these estimates together we thus have \beqno{g1a}\sup_{\eg>0}|\mg_1*\fg_\eg| \le C
\Psi_2^\frac{1}{2}\Psi_3^\frac{1}{s'}.\eeq

By H\"older's inequality, we get
$$\iidx{\Og}{\sup_{\eg>0}|\mg_1*\fg_\eg|} \le \left(\iidx{\Og_R}{\Psi_2}\right)^\frac1{2}\left(\iidx{\Og_R}{\Psi_3^\frac{2}{s'}}\right)^\frac{1}{2}. $$

We have
$$\left(\iidx{\Og_R}{\Psi_2}\right)^\frac{1}{2} = \left\|M\left(|H|^{2(p-1)}|\dsg H|^{2}\right)\right\|_{L^1(\Og_R)}^\frac12.$$

As $s'=\frac{2N}{N+2\sg}<2$, we can also use \mref{HL} with the Lebesgue measure $wdx=dx$ to get from the definitions of $\mI_1, h$, $$ \left(\iidx{\Og_R}{\Psi_3^\frac{2}{s'}}\right)^\frac{1}{2}\le C\left( \iidx{\Og_R}{[|H|^{p}\dsg u]^{2}} \right)^{1/2} = C\mI_1^{1/2} .$$

Therefore,  \beqno{g1est}\iidx{\Og_R}{\sup_\eg|\mg_1*\fg_\eg|} \le C \mI_1^\frac12\left\|M\left(|H|^{2(p-1)}|\dsg H|^{2}\right)\right\|_{L^1(\Og_R)}^\frac12.\eeq

The term $\sup_\eg |\mg_2 *\fg_\eg|$ is handled similarly (again, replacing $h$ by its average $\mitx{B_{\eg}}{h}$ in $\Psi_3$). We then obtain an estimate like \mref{g1est} to have $$\iidx{\Og_R}{\sup_{\eg>0}|(\mg_1+\mg_2)*\fg_\eg|} \le C \mI_1^\frac12\left\|M\left(|H|^{2(p-1)}|\dsg H|^{2}\right)\right\|_{L^1(\Og_R)}^\frac12.$$
This is \mref{newg1estzz} and the lemma is proved. \eproof

Of course, it is easy to remove the use of maximal function $M$ in \reflemm{g12lemM} and have
\blemm{g12lemr} Suppose that  $2<r$, $r'> s'$.  If we define
$$\mI_{1,r}=\iidx{\RR^N}{|H|^{(p-1)r'}(x)|\dsg u|^{r'}(x)\og^{r'}}\quad \mI_{2,r}=\iidx{\RR^N}{|H|^{(p-1)r}(x)|\dsg H|^{r}(x)\og^r}  ,$$
Then \reflemm{g12lemM}  holds in the form $$\iidx{\Og_R}{\sup_{\eg>0}|(\mg_1+\mg_2)*\fg_\eg|\og^2} \le \cps(N,\sg) \mI_{1,r}^\frac{1}{r'}\mI_{2,r}^\frac{1}{r}.$$ \elemm
\bproof Indeed, for such $r$, we can use H\"older's inequality and \mref{g1a} to have ($\Psi_2, \Psi_3$ are defined as in the proof of \reflemm{g12lemM})
$$\iidx{\Og}{\sup_{\eg>0}|\mg_1*\fg_\eg|} \le \left( \iidx{\Og_R}{\Psi_2^\frac{r}{2}}\right)^\frac{1}{r}\left( \iidx{\Og_R}{\Psi_3^\frac{r'}{s'}}\right)^\frac{1}{r'}. $$
As $\frac{r}{2}>1$ we can use \mref{HL} for the first integral to estimated by $\mI_{2,r}$ defined here. Finally, because  $r'/s'> 1$, we apply \mref{HL} again for the second one on the right hand side to obtain the claim. \eproof

We also need to estimate $\|\mg_1+\mg_2\|_{L^1}$ but this can be done similarly. Indeed, considering the map $T:f\to (\mg_1+\mg_2)*f$ on $L^1$, we easily see that $\lim_{\eg\to0}\|T\fg_\eg\|_{L^1}\ge \|\mg_1+\mg_2\|_{L^1}$. The estimate for $\|T\fg_\eg\|_{L^1}$ can be obtained by the previous lemmas.

Finally, to completely remove  $r$ from $\mI_{2,r}$ we can follow \mref{Idefz} and recall the definition of  $\mI_2$ in \mref{Idefz}
$$
\mI_2:=\iidx{\RR^N}{|H|^{2p-2}|\dsg H|^2\og^2}.$$

\blemm{mainlocWBMOlem} There is a contants $C$ such that
\beqno{GNglobalestz00}\mI_1\le C(N,p,\sg)\|u\|_{BMO(\Og_R)}^2 \mI_2+ E_\Ofam(u)
\iidx{\Og}{|H|^{2p}}+C\iidx{\Og}{|H|^{2p}\dsg \og|^2}.\eeq
\elemm

The proof will make use of different integration by parts   via \reflemm{lapintlem}.

\bproof First, we establish that for $0<\bg<\sg/2$ and   $h=|H|^{p-1}H$ \beqno{GNglobalestzab}\mI_1\le C\|\dgeng{\bg} u\|_{BMO(\Og_R)}^2 \iidx{\RR^N}{|\dsg h|^2\og^2}+E_\Ofam(u)
\iidx{\Og}{|h|^{2}}+C\iidx{\Og}{|H|^{2p}\dsg \og|^2}.\eeq

Recall that $\mg_i=\dgeng{\sg} (U_i)$. For any $y\in\RR^N$ and $\eg>0$ (note that $U_i$ are the same as before)
$$U_1= h\dsg u\left(h-\mitx{B_{\eg}(y)}{h}\right),\;U_2= (h-\mitx{B_{\eg}(y)}{h})\dsg u\mitx{B_{\eg}}{h},$$
$$U_3=\dsg u\left(\mitx{B_{\eg}(y)}{h}\right)^2.$$

Again as in \mref{keyGNwsp} (again, via integration by parts and dropping conjugates for simplicity), if $\sg=\ag+\bg$ then we have \beqno{I1est}\barrl{\mI_1=\iidx{\RR^N}{\dgeng{\ag}(|H|^{2p}\dsg u \og^2) \dgeng{\bg} u}\le }{1.5cm} &\left|\iidx{\RR^N}{(\mg_1+\mg_2) \dgeng{\bg} u\og^2}\right|+\left|\iidx{\RR^N}{\mg_3 \dgeng{\bg} u\og^2}\right|+
\\&
\iidx{\RR^N}{|H|^{2p}|\dsg u||\og||\dsg\og|}. \earr\eeq

To estimate the first  term on the right, we argue similarly as before. 
First of all, by the Fefferman-Stein theorem again ($\mbox{supp}\og\subset\Og_R$), $$\left|\iidx{\RR^N}{(\mg_1+\mg_2) \dgeng{\bg} u\og^2}\right|\le \cfs(N)\|\dgeng{\bg} u\|_{BMO{(\Og_R)}}\|(\mg_1+\mg_2)|\|_{\ccH^1(\Og_R)}.$$ 

For any $y\in\Og$, we use integration by parts (via Parseval's theorem again), ignoring $\og^2$ for simplicity  because the integrals involving $\og, \dsg \og$ in \mref{keyGNwsp} can be treated the same way as we estimate  $\mg_3$ in \reflemm{g3lem} later (after an integration by parts in $\dgeng{\ag}$). 

Now, the property of $\fg_\eg$ in \mref{fgprop} is changed accordingly as we replace $\sg$ in the proof of \reflemm{g12lemM} by $\ag$. Together with H\"older's inequality for any $s>1$, we have the following
$$\barr{lll}|\mg_1*\fg_\eg(y)|&=&\left|\iidx{B_\eg(y)}{\dgeng{\ag} \fg_\eg(x-y)(h-h_{B_\eg(y)})h\dsg u  }\right|\\&\le&\frac{C_1}{\eg^{\ag}}
\left|\mitx{B_\eg(y)}{|h-h_{B_\eg(y)}||h\dsg u|}\right|
\\&\le&\frac{C_1}{\eg^{\ag}}
\left(\mitx{B_\eg(y)}{|h-h_{B_\eg(y)}|^s}\right)^\frac1s 
\left(\mitx{B_\eg(y)}{|h\dsg u|^{s'}}\right)^\frac1{s'}.
\earr$$

Apply the Sobolev inequality to the first integral on the right with $\sg, s$ being replaced by $\ag$, $\frac{2N}{N-2\ag}$  respectively. We still have $s_*=2$ (this is essential)  and $|\dgeng{\ag} h| \lesssim c(p)|H|^{p-1}|\dgeng{\ag}H|$, power rule inequality in \refrem{prodineq}.

Applying \reflemm{lapintlem} twice, integration by parts, for any smooth function $f$ we have 
$$\iidx{\RR^N}{\myprod{ f, \overline{ \dgeng{\ag+\bg} u} }}=\iidx{\RR^N}{\myprod{\dgeng{\ag} f, \overline{\dgeng{\bg} u} }}=\iidx{\RR^N}{\myprod{ f, \overline{ \dgeng{\ag}\left( \dgeng{\bg} u\right) } }}.$$

Thus, $\dgeng{\sg-\ag} (\dgeng{\ag} h)=\dsg h$. Suppose that $\dsg h\in L^2(\RR^N)$ then $\dgeng{\sg-\ag}f \in L^2(\RR^N)$ with $f=\dgeng{\ag} h$. The Sobolev inequality for $f$ in $W^{\sg-\ag,2}(\RR^N)$ then implies
$$\left(\iidx{\RR^N}{|\dgeng{\ag} h|^r}\right)^\frac 2r \le C \iidx{\RR^N}{|\dgeng{\sg-\ag} f|^2}=C \iidx{\RR^N}{|\dsg h|^2}, \quad r=\frac{2N}{N-2(\sg-\ag)}>2. $$

Thus we can find $r>2$ such that for $\ag<\sg$ (note that we are ignoring $\og$) $$\left(\iidx{\RR^N}{|\dgeng{\ag} h|^r} \right)^\frac 2r \le C\iidx{\RR^N}{|\dsg h|^2}. $$ Hence, we can use \reflemm{g12lemr}, with $\sg$ being replaced by $\ag$, as  $r'=\frac{2N}{N+2(\sg-\ag)}> s'=\frac{2N}{N+2\ag}$ also, if $2\ag>\sg$ or $\bg<\sg/2$.  By the above inequality and because $r'<2$, we can  replace $\mI_{1,r}, \mI_{2,r}$ respectively  by $\mI_1,\mI_2$ and we obtain \mref{GNglobalestzab} with $0< \bg<\sg/2$.  We can let $\bg\to0^+$ in \mref{GNglobalestzab} to obtain \mref{GNglobalestz00}. In fact, we see that if $\Og$ is an extended bounded domain, then $\lim_{\bg\to0^+}\|\dgeng{\bg}u\|_{BMO(\Og)}\le \|u\|_{BMO(\Og)}$. Indeed, as $\dgeng{\bg}u\to u$ (again using the Fourier transformation definition, see \cite[Proposition 4.4]{NPV},  and the dominated convergence theorem as $u\in W^{\sg,2}(\Og)$), we see in the following that $\lim_{\bg\to0^+}(\dgeng{\bg}u)_{B_R}=(u)_{B_R}$ for all balls $B_R$. This is not a simple fact. First of all, we easily see that for all $r>0$
$$\barr{lll}\left| \idz{B_r(x)}{\frac{u(x)-u(y)}{|x-y|^{N+\bg}}}{dy} \right| &=& \left| \idz{B_r(x)}{\frac{1}{|x-y|^{N}}\frac{u(x)-u(y)}{|x-y|^{(\bg+\dg)-\dg}}}{dy} \right|\\
&\le&  \idz{B_r(x)}{\left( \frac{|u(x)-u(y)|^2}{|x-y|^{N+2\bg+2\dg}}+\frac{1}{|x-y|^{N-2\dg}}\right)}{dy}. \earr    $$

Obviously, the last integral is a $L^1(\Og)$ function (in $x$) if $r< 1, \dg>0$ and all $\bg$ such that  $0<\bg+\dg<\sg$, as  $u\in W^{\sg,2}(\Og)$. In addition,
$$\left|\idz{\Og\setminus B_r(x)}{\frac{u(x)-u(y)}{|x-y|^{N+\bg}}}{dy}\right| \le \frac{1}{r^{N+\bg}}\idz{\Og\setminus B_r(x)}{|u(x)-u(y)|}{dy}, $$ which is  in $L^1(\Og)$
because $ u\in L^2(\Og)$, and using Fubini's theorem. Thus, we can use the dominated convergence theorem as we just proved that $|\dgeng{\bg}u|\le g$ for some $g\in L^1(\Og)$. We also  note \cite{BBM, MazS}, where it showed that the Gagliardo norm $\|u\|_{W^{\bg,p}(\Og)}$ (not $\|\dgeng{\bg}u\|_{L^1(\Og)}$ when $u\in W^{\bg,2}(\Og)$) may blow up when $\bg\to 0^+$ if the scaling factor is not appropriate.

We now  can use the definition of the BMO norm $\|u\|_{BMO(\Og)}=\sup_{B_R}\frac{1}{|B_R|}\iidx{B_R}{|u-(u)_{B_R}|}$ to see that $\lim_{\bg\to0^+}\|\dgeng{\bg}u\|_{BMO(\Og)}\le \|u\|_{BMO(\Og)}$, this is elementary.

 We will estimate the last two terms in \mref{I1est} later when we consider $\mathfrak{J}_2$ to completely obtain \mref{GNglobalestz00}.  \eproof

Concerning the last term $\mathfrak{J}_2$ involving $\mg_3$ in \mref{keyGNwsp}, we have the following lemma. Notice that we do not have to estimate ${\cal H}^1$ norms although there is $\eg$ in the definition of $\mg_3$.

\blemm{g3lem} For any $\eg, \thg>0$ we have
$$\mathfrak{J}_2\le E_\Ofam(u)\iidx{\Og}{h^2}+\thg\iidx{\Og}{|H|^{2p}|\dsg u|^2\og^2}+C(\thg)\iidx{\Og}{|H|^{2p}|\dsg \og|^2}.$$
where (see \mref{E*def} and the definition of the family $\Ofam$)  $$E_\Ofam(u)=\sup_{B^i\in\Ofam}|B^i|^{-1}\iidx{B^i}{[|\dsg u|^2+|u|^2\LLg_i^2]}.$$
\elemm

\bproof From  the definition of $\Ofam$, we first have by integration by parts on each $B^i$ and then Young's inequality that
$$\barr{lll}\left|\iidx{\Og}{\mg_3 u\og^2}\right|&\le& \sum_i\left| \iidx{B^i}{|h_{B^i}|^2\dsg u u\og_i^2}\right|\\
&\le& \sum_i\iidx{B^i}{|h_{B^i}|^2[|\dsg u|^2\og_i^2+|u|^2|\LLg_i|^2]}.\earr$$

The right hand side is estimated by $$\sum_i|h_{B^i}|^2\iidx{B^i}{[|\dsg u|^2\og_i^2+|u|^2\LLg_i^2]}\le \sup_i\iidx{B^i}{[|\dsg u|^2\og_i^2+|u|^2\LLg_i^2]}\sum_i|h_{B^i}|^2.$$
Using the inequality $|h_{B^i}|^2\le (h^2)_{B^i}$ and the finite intersection property of $B^i$'s we bound the last term by
\beqno{g3esta} \sup_i|B^i|^{-1}\iidx{B^i}{[|\dsg u|^2\og_i^2+|u|^2\LLg_i^2]}\iidx{\Og}{h^2}.\eeq

Regarding the second term is $\mathfrak{J}_2$, we use Young's inequality
$$\iidx{\Og}{|H|^{2p}|\dsg u||\og||\dsg\og|}\le \iidx{\Og}{|H|^{2p}[\thg|\dsg u|^2\og^2+C(\thg)|\dsg \og|^2]}.$$

Combining the above estimates, we obtain the lemma. \eproof

\brem{J2rem} The reason we did not use a partition in estimating the second term in $\mathfrak{J}_2$ is that we will obtain a better estimate because the number of the balls $B^i$'s would be very large. Thus, leaving $\og$ as it is would be a better choice. Also, we separate $E_\Ofam$ from the second term in $\mathfrak{J}_2$ to freely choose $\Ofam$ from $\og$ later on. It is important to note that the factor $E_\Ofam$ is multiplied by the  integral of $|h|^2=|H|^{2p}$ over $\Og$.  These facts do not effect the local qualitative analysis when the main goal is to absorb $\mI_1$ into $\mI_2$ and this requires only that $\|u\|_{BMO}$ is small.
\erem

\brem{EOfam} In the proof of \reflemm{g3lem}, we abused the product rule for $\dsg$ which is quite complicated (see \cite{zua}) but it is obvious when $\sg=1$, to be more precise $E_\Ofam \sim |B^i|^{-1}\|\dsg (u\og_i)\|_{L^2(\RR^N)}$ (we can assume that $u$ is bounded because this term will not appear eventually). However, we eventually prove that $E_{\Ofam}\to0$ by choosing $\Ofam$ properly (via the partition $\og_i$'s whose $\dgeng{\sg}\og_i$'s are independent of $\dgeng{\sg}\og$) and obtain a much better version of \reftheo{GNlobalz} without $E_\Ofam$ in \refcoro{locwGNcoro}. Without being concerned with these technicalities to present the main ideas, we leave the proof as it is and the details for interested readers.
\erem

Finally, we are ready to have

{\bf Proof of \reftheo{GNlobalz}:}
The last term in \mref{keyGNwsp} is treated  in \reflemm{g3lem}. Putting the estimates in \reflemm{mainlocWBMOlem} and \reflemm{g3lem} together and choose $\thg$ sufficiently small, we obtain 
\reftheo{GNlobalz}. \eproof

\subsection{A variances of the weak inequality in \reftheo{GNlobalz} by a choice of $\Ofam$} \label{varWGN}

{\bf An estimate for $|\dsg\fg|$ of a Lipschitz function $\fg$:} Due to the presence of $\LLg_i$'s in $E_\Ofam$, we need some elementary estimate for $|\dsg\fg|$ of a Lipschitz function $\fg$ (or $\og_i$) in order to improve the choice of $\Ofam$ later on. Let $\LLg,k>0$ and $\fg$ be the function
$$\fg=\LLg\chi_{B_{R+k^{-1}}}.$$

We have $\dsg \fg=\mathfrak{F}^{-1}(|\zeta|^\sg\mathfrak{F}\fg)$.
Making a change of vars $\zeta=k\rg$ and recalling the property that $k^N\mathfrak{F}(\fg)(k\rg)=\mathfrak{F}(\mbox{dil}_{k^{-1}}\fg)(\rg)$ ($\mbox{dil}_\tau$ is the dilation operator) we get
$$ \barr{lll} \dsg \fg(x) &=& \int_{\RR^N} e^{i2\pi x\zeta}|\zeta|^\sg \mathfrak{F}(\fg)(\zeta) d\zeta =  k^\sg\int_{\RR^N} e^{i2\pi x\cdot k\rg}|\rg|^\sg k^N\mathfrak{F}(\fg)(k\rg) d\rg\\
&=& k^\sg\dspl{\int}_{\RR^N} e^{i2\pi x\cdot k\rg}|\rg|^\sg \mathfrak{F}(\mbox{dil}_{k^{-1}}\fg)(k\rg) d\rg.\earr
$$
Clearly, $\mbox{dil}_{k^{-1}}\fg = \LLg\chi_{B_{1+kR}}$, and $\mbox{dil}_{k^{-1}}\fg \to \LLg\chi_{B_{1}}$ as $k\to0$, and $\mbox{dil}_{k^{-1}}\fg \to \LLg$ as $k\to\infty$. They have compact supports for $0<k<\infty$ and uniformly bounded. So that $\mbox{dil}_{k^{-1}}\fg$ is  a Schwartz function (by approximation) and so is the last integral.
Hence, $|\dsg \fg|\le c\LLg k^\sg$. 

If $\LLg\sim k^{-1}$ and $k\in(0,C]$, then  there is a constant $c$ such that
\beqno{dsgog}|\dsg \fg|\le ck^{\sg-1} . \eeq

Of course, the above is not true if  $k\to\infty$ because the above calculation does not apply as $\mbox{dil}_{k^{-1}}\fg$ is no longer a Schwartz function and uniformly bounded (and $c$ does not exist).

We are going to use \mref{dsgog} to show that this is the case if we assume $u, \dsg u\in L^2(\RR^N)$ and $H\in L^{2p}(\Og)$ (this is not assumed  in \reftheo{GNlobalz} to obtain  the local weak inequality and we are going to let $\Og, B^i\to\RR^N$) by choosing $B^i$'s (in the family of covering) $\Ofam$  and $\og$ appropriately.

Recall that (see \mref{E*def} and the definition of the family $\Ofam$)  $$E_\Ofam(u)=\sup_{B^i\in\Ofam}|B^i|^{-1}\iidx{B^i}{[|\dsg u|^2+|u|^2\LLg_i^2]}.$$

We now look carefully at the following two terms in $\mathfrak{J}_2$ and the choices of $\Ofam, \og$ $$E_\Ofam(u)\iidx{\Og}{|H|^{2p}} \mbox{ and } \iidx{\Og}{|H|^{2p}|\dsg \og|^2}.$$

Assume that $u\in BMO(B_R)$, $k<1$, and we define $\og$ to be a Lipschitz function 
$$\og(x)=k^{-1}\left\{\barr{ll}1&x\in B_R,\\ 0&x\in \RR^N\setminus B_{R+k^{-1}} \earr\right.$$
By \reftheo{GNlobalz}, (multiplying both sides of the inequality by $k^2$) we have $$\iidx{B_{R}}{|H|^{2p}|\dsg u|^2\og^2}\le C\|u\|_{BMO(B_R)}^2\iidx{B_{R+k^{-1}}}{|H|^{2p-2}|\dsg H|^2\og^2} +k^2[\ldots],$$ where $k^2[\ldots]\to0$ when $k\to0$, see \mref{E*def} and the definition of the family $\Ofam$ when $B^i$'s are balls with radius $R+k^{-1}$ (so that $\LLg_i^2\le ck^{2\sg-2}$), and $k^2|\dsg \og|^2 \le ck^{2\sg}$).

From the proof of \reflemm{g3lem}, we  have by integration by parts ($\LLg_i=|\dsg\og_i|$)
$$E_\Ofam \le c(\|u\|_{L^\infty(\Og)}i)\sup_{B^i\in\Ofam}|B^i|^{-1}\left[ \iidx{B^i}{[|\dsg u|^2\og_i^2+|u|^2|\LLg_i|^2]}+\|\dsg(u\og_i)\|_{L^2(B_i)}^2\right].$$

As $\dsg u, u\in L^2(\RR^N)$ and $\dsg (u\og_i)\in L^2(\RR^N)$ because $\og_i$ is Lipchitz (see \cite{NPV}), it is easy to use the Fourier transform definition of $\dsg$ to prove that $\lim_{k\to0}|B_i|^{-1}\|\dsg(u\og_i)\|_{L^2(B_i)}^2=0$ (as $\|\dsg(u\og_i)\|_{L^2(B_i)}^2\le (1+\LLg_i^2)\|u\|^2_{W^{\sg,2}(\RR^N)}$, we omit the simple argument and this implies $E_\Ofam\to0$).

Now, if $H, u$ can be extended to $\RR^N$ and we have $$u, \dsg u\in L^2(\RR^N),\; |H|^{2p} \in L^1(\RR^N),$$ we need only $u\in BMO(B_R)$, then we can let $k\to0$ and then $\mbox{supp}\og\to\Og(\RR^N), R\to\infty$ in the above to see that (as $E_\Ofam \sim |B^i|^{-1} \|\dsg (u\og_i)\|_{L^2(\RR^N)}^2\to0$ and $|\dsg \og_i|\to0$)
$$\iidx{\RR^N}{|H|^{2p}|\dsg u|^2}\le C\|u\|_{BMO(B_R)}^2\iidx{\RR^N}{|H|^{2p-2}|\dsg H|^2}.$$

We thus have a local and global result as follows ($H^{2p}\to H^{2p}\og^2$ )
\bcoro{locglobWGNcoro} Assume that $R>0$ and measurable $H,u$ on $B_R$ are such that $u\in BMO(B_R)$ and $H, u$ can be extended to $\RR^N$ such that $u, \dsg u\in L^2(\RR^N),\; |H|^{2p} \in L^1(\RR^N)$. Then,
\beqno{GNglobalestzcoro}\iidx{\RR^N}{|H|^{2p}|\dsg u|^2}\le C(N,p,\sg)\|u\|_{BMO(B_R)}^2\iidx{\RR^N}{|H|^{2p-2}|\dsg H|^2}.\eeq
The above integrals can be both infinite. If $|H|^{2p-2}|\dsg u|^2 \in L^1(\RR^N)$ then they are both finite.
\ecoro

The above argument does not apply if $k\to\infty$ so that we can not replace the integral on the right hand side  over $\RR^N$ by the one over $B_R$ (even $H\in W^{\sg,2}(B_R)\cap L^{2p-2}(B_R)$). Of course, the left hand side integral can be reduced to that over $B_R$. A reason is that if  $H\in W^{\sg,2}(B_R)\cap L^{2p-2}(B_R)$ then we can extend $H$ to $\RR^N$ and  have a constant $C$ such that
$$\iidx{\RR^N}{|H|^{2p-2}}\le C\iidx{B_R}{|H|^{2p-2}},\; \iidx{\RR^N}{|\dsg H|^{2}}\le C\iidx{B_R}{|\dsg H|^{2}}$$ but we can not find such a constant $C$ and
$$\iidx{\RR^N}{|H|^{2p-2}|\dsg H|^2}\not\le C\iidx{B_R}{|H|^{2p-2}|\dsg H|^2}.$$

If $u\in BMO(B_R)$ and $u, |H|^{p-1}H\in W^{\sg,2}_0(B_R)$ then by extending $H$ to be zero outside $B_R$ we obtain
$$\iidx{B_R}{|H|^{2p}|\dsg u|^2}\le C\|u\|_{BMO(B_R)}^2\iidx{B_R}{|H|^{2p-2}|\dsg H|^2}.$$
In fact, this is also true if $|H|^{p-1}H\in W^{\sg,2}(B_R)$ as  the last integral is in fact that of $|\dsg (|H|^{p-1}H)|^2$ if we inspect the proof of  \reflemm{g12lemM} where we estimate $\mg_1$. Thus, $\og$ and $\mathfrak{J}_2$ in our local weak inequality of \reftheo{GNlobalz} can be removed, keeping only $\mI_1, \mI_2$, if $\Og$ is an extended domain (as $B_R$). However, due to the common nature of the local analysis of the regularity of weak solutions to PDE's, one has $\og$ in the consideration and \reftheo{GNlobalz} is sufficient. 

The fact that $E_\Ofam\to0$, when $k\to0$, greatly improves \reftheo{GNlobalz} so that we state it as follows.

\bcoro{locwGNcoro} Assume that $\bg\ge0$, $\sg_H\ge \sg_u-\bg$. If $\dgeng{\bg}u\in BMO(\Og)$, $u \in W^{\sg_u,2}(\Og)$, $h=|H|^{p-1}H\in W^{\sg_H,2}(\Og)$, and $\Og$ is an extended domain, then $$\iidx{\Og}{|H|^{2p}|\dgeng{\sg_u} u|^2}\le C(N,p,\sg_u,\Og)\|\dgeng{\bg}u\|_{BMO(\Og)}^2\||H|^{p-1}\dgeng{\sg_H}H\|_{L^{2}(\Og)}^2.$$
The above inequality also holds when $\Og=\RR^N$ and, of course, $C(N,p,\sg,\RR^N)=C(N,p,\sg)$. 
\ecoro

One should note that,  according to the Fourier transform definition, the last factor is zero iff $H\equiv0$.
 Also, the constant $C(N,p,\sg,\Og)$ depends on $\Og$ due to the constant $C(N,\Og)$ in the usual extension map $\|g\|_{W^{\sg,2}(\RR^N)}\le C(N,\Og)\|g\|_{W^{\sg,2}(\Og)}$.
 
If $|H|^{p-1}H, |H|^{p-1}\dsg H\in L^2(\Og)$, then $|H|^{p-1}H\in W^{\sg,2}(\Og)$. The above corollary also  yields
$$\iidx{\Og}{|H|^{2p}|\dsg u|^2}\le C(N,p,\sg,\Og) \|u\|_{BMO(\Og)}^2\iidx{\Og}{|H|^{2p-2}|\dsg H|^2}.$$

The integral on the left hand side is finite. This assertion is not trivial under the assumptions that $u\in BMO(\Og)$ and $u, |H|^{p-1}H\in W^{\sg,2}(\Og)$ (even when $H=u$). When $p=1$, we also have
$$\mbox{$u\in BMO(\Og)$, and $u, H\in W^{\sg,2}(\Og)$}\Rightarrow\iidx{\Og}{|H|^{2}|\dsg u|^2}\le C(N,\sg,\Og)\|u\|_{BMO(\Og)}^2\|H\|_{W^{\sg,2}(\Og)}^2.$$

This leads to the special case when $H=u$. We can easily iterate \refcoro{locwGNcoro} to assert that:
$$\mbox{If  $H\in W^{\sg,2}(\Og)\cap BMO(\Og)$,  then for all $p\ge1$ }\iidx{\Og}{|H|^{2p}|\dsg H|^2}<\infty.$$
Also, the above implies that if   $H\in W^{\sg,2}(\Og)\cap BMO(\Og)$,  then for all $p\ge1$ $|H|^{p-1}H$ can be extended to $W^{\sg,2}(\RR^N)$. 

As a byproduct of \refcoro{locwGNcoro}, we also have the following theoretic result. If $\Og$ is a bounded extension domain and $u\in BMO(\Og)$ and $u,H\in W^{\sg,2}(\Og)$, then by H\"older's inequality for $p>1$
$$\iidx{\Og}{|u|^2|\dsg H|^2}\le \||u|^2\|_{L^{p'}(\Og)} \||\dsg H|^2\|_{L^{p}(\Og)}.$$

By the continuity of the map $p\to \|g\|_{L^p(\Og)}$ for a given $g\in L^p(\Og)$ and $u\in L^q(\Og)$ for all $q>1$, we conclude that
$$\iidx{\Og}{|u|^2|\dsg H|^2}\le C(\|u\|_{BMO(\Og)}, \|\dsg H\|_{L^{2}(\Og)})<\infty.$$
Together with \refcoro{locwGNcoro}, if $\sg=1$, then we now see that $uH\in W^{1,2}(\Og)$. Thus,

\bcoro{Walgeb} Assume that $\Og$ is a bounded extension domain and that $u,H\in W^{1,2}(\Og)$. If either $u\in BMO(\Og)$ or $H\in BMO(\Og)$, then $uH\in W^{1,2}(\Og)$. 
\ecoro

By iterating \mref{GNglobalestzcoro}  in $p$, if $u\in BMO$ and $|H|^p, \dsg H\in L^2(\RR^N)$, then $|H|^p\dsg u\in L^2(\RR^N)$ for all $p\ge1$. Hence, (with $H=u$) we have for $H\in BMO(\RR^N)\cap W^{\sg,2}(\RR^N)\cap L^2(\RR^N)$
\beqno{GNglobalestzcoroHisu}\iidx{\RR^N}{|H|^{2p}|\dsg H|^2}\le C(N,p,\sg)\|H\|_{BMO(\RR^N)}^2 \iidx{\RR^N}{|H|^{2p-2}|\dsg H|^2}. \eeq
This trivial if $H\in L^\infty(\RR^N)\cap W^{\sg,2}(\RR^N)$.



\bibliographystyle{plain}

\end{document}